\documentclass[12pt]{amsart}
\usepackage{a4wide,enumitem,amsmath,amsfonts,amssymb,amsthm,graphics,amsrefs,amscd}
\usepackage{hyperref}
\usepackage{amsthm,amsxtra}

\makeindex

\newcommand{\Z}{\mathbb{Z}}

\newcommand{\C}{\mathbb{C}}
\newcommand{\Q}{\mathbb{Q}}

\renewcommand{\Re}{\operatorname{Re}}

\newtheorem*{theorem}{Theorem}
\newtheorem*{lemma}{Lemma}

\newtheorem*{proposition}{Proposition}

\newtheorem*{example}{Example}
\newtheorem*{definition}{Definition}
\newtheorem*{remark}{Remark}
\newtheorem*{claim}{Claim}
\newtheorem*{principle}{Principle}

\numberwithin{equation}{section}

\def\ms{{\medskip}}
\def\sms{{\smallskip}}

\def\del{{\delta}}
\def\al{{\alpha}}
\def\lam{{\lambda}}
\def\sig{{\sigma}}
\def\om{{\omega}}

\def\Gam{{\Gamma}}
\def\Lam{{\Lambda}}
\def\eps{\epsilon}

\def\Sig{{\Sigma}}

\def\chG{{\mathcal{G}\spcheck}}
\def\chT{{\mathcal{T}\spcheck}}
\def\chR{{\mathcal{R}\spcheck}}
\def\chS{{\mathcal{S}\spcheck}}
\def\chH{{\mathcal{H}\spcheck}}

\def\cGD{{\mathcal{G}^D}}

\def\chSig{{\Sig\spcheck}}

\def\cA{{\mathcal{A}}}

\def\cB{{\mathcal{B}}}

\def\cF{{\mathcal{F}}}
\def\cR{{\mathcal{R}}}
\def\cP{{\mathcal{P}}}
\def\cM{{\mathcal{M}}}
\def\cH{{\mathcal{H}}}
\def\rC{{\mathcal{C}}}

\def\GD{{\mathcal{G}^D}}

\def\bW{{\mathbf{W}}}
\def\bG{{\mathbf{G}}}
\def\bX{{\mathbf{X}}}
\def\bL{{\mathbf{L}}}
\def\bDel{{\mathbf{\Delta}}}
\def\bA{{\mathbf{A}}}
\def\bT{{\mathbf{T}}}
\def\bH{{\mathbf{H}}}
\def\bR{{\mathbf{R}}}

\def\bS{{\mathbf{S}}}
\def\bR{{\mathbf{R}}}

\def\GalF{{{\rm Gal}(F)}}
\def\WF{{W_F}}
\def\tWF{{\tilde{W}_F}}

\def\LANG{{LAN(G)}}
\def\cLANG{{\mathcal{LAN}(G)}}

\def\IrrH{{\rm{Irr}(H)}}

\def\Fbar{{\bar{F}}}

\def\MG{{\mathcal{M}(G)}}
\def\MT{{\mathcal{M}(T)}}
\def\MGe{{\mathcal{M}(G,\eps)}}

\def\MHe{{\mathcal{M}(H,\eps)}}

\def\tbG{{\tilde{\mathbf{G}}}}

\def\tT{{\tilde{T}}}
\def\tG{{\tilde{G}}}

\def\BG{{\mathcal{GR}(G)}}

\begin{document}

\title[Hidden sign]{Hidden sign in the Langlands correspondence}
\author{Joseph Bernstein}
\address{School of Mathematical Sciences, Tel Aviv University, Tel Aviv 6997801, Israel}
\email{bernstei@tauex.tau.ac.il}

\begin{abstract}

  In this note I describe some modification of the Langlands correspondence and explain why the modified correspondence is more natural than the original one.
  I will also discuss the significance of this modification to the theory of automorphic $L$-functions.

\end{abstract}

\maketitle


\section{Introduction}

\subsection{The Langlands correspondence}

In 1967, motivated by his results on Eisenstein series, Langlands introduced the notion of the \textit{dual group}. 
Namely, starting with a~reductive group $G$ 
over some field $F$, he constructed the dual group $\chG$ over $\C$. He suggested that many problems 
about representations of groups over local fields, and about automorphic representations of  adelic groups,
should have natural interpretations in terms of the Galois group of $F$, and of the dual group $\chG$.

Langlands  formulated several conjectures that described the correspondence between representations 
of $G$ and Galois objects constructed in terms of the dual group $\chG$ (the \textit{Langlands correspondence}).
This discovery transformed the area of representation theory and automorphic forms since it 
gave some arithmetic flavor 
to the problems of representation theory of these groups.
During the last 50 years, this paradigm (under the name ``the Langlands program") is one of the main lines of 
research in the Representation Theory and the Theory of Automorphic Representations.   

   In this note I will discuss only the $p$-adic case of the Local Langlands correspondence (LLC), though 
   the conclusions are clearly applicable
   to representations of real groups, to automorphic representations, and to the Geometric Langlands Program.
   
   \ms
   
   In what follows, we fix a~$p$-adic field $F$ of characteristic $0$. We denote by $\WF$ the \textbf{Weil group}
    of the field $F$ -- the subgroup of the absolute Galois group $\GalF := {\rm Gal}(\bar{F} / F)$ that in some sense better 
   captures the arithmetic
   structure of $F$. By definition, $\WF$ consists of elements of the group $\GalF$ that act on the residue field 
   by integer powers of the  Frobenius automorphism.
   
    We fix a~connected reductive group $\bG$ over $F$,  and consider the topological group ${G = \bG(F)}$ of its $F$-points.
    Our main object of study is the category $\MG$ of \textbf{smooth} representations of the group $G$ (we start 
    with  representations in vector spaces over complex numbers, but later 
    switch to slightly more general case). 
    
    One of the most important problems here is to give a~good description of the set ${\rm Irr}\ G$ of irreducible 
    representations of $G$, i.e., the set of isomorphism classes of simple objects in $\MG$.
    
   In this note I will discuss only the case of a~split group $G$.
    Langlands constructed the complex reductive group $\chG$ dual to $\bG$ and proposed to classify 
    irreducible representations of $G$ in terms of $\chG$.  
    
    Namely, he considered the set of Langlands parameters 
    $\LANG$ that consists of continuous morphisms $\phi: \WF \to \chG$ modulo the adjoint action of the group $\chG$.
    The main claim of the Local Langlands Correspondence (LLC) is that there should be some natural correspondence
     between the sets  ${\rm Irr}\ G$ and $\LANG$.
     
     \sms
     
     Later on, the correspondence of this type has been constructed in many cases. It has been also realized that, in
     order to have such a~natural correspondence, one has to modify both sides of it.
     
     On the representation side, ${\rm Irr}\ G$ should be replaced by the set of isomorphism classes of stacky $G$-modules 
     (see \cite{B}). In more explicit terms, this means that instead of the set ${\rm Irr}\ G$ we should
      consider the disjoint union of the sets    ${\rm Irr}\ G_i$ for all pure inner forms of the group $G$.
     
     On the Galois side, the quotient $\LANG = {\rm Mor}(\WF, \chG) / \chG$
     should be considered not as a~set, but as an algebraic stack $\cLANG$ (let us call it the \textbf{Langlands stack}).
     The set of Langlands parameters should be extracted from this stack by a~standard procedure.
     
     \sms
     
     In this note I will concentrate on the Galois side of the LLC. I will argue that, in fact, there is some hidden sign
      in the Langlands Correspondence. This implies that  in order to have a~more natural correspondence, one has to modify the Langlands stack $\cLANG$ described above.
      
      \subsection{ }
      
      Let me describe some general semantic difficulty in analyzing this type of correspondences.
      Suppose we have two constructions that produce two sets $S$ and $T$, and we expect that there is 
      a~``natural" correspondence between these two sets. For example, we would like to construct a~map
      $\nu: S \to T$ which is a~bijection, or is very close to a~bijection.
      
      Suppose that after some work I produced a~nice map $\nu: S \to T$. However, suppose that my colleague, using different 
      technique, has produced  another map $\mu: S \to T$. What should I think about this?
      
 The first idea is that, in fact, these two maps should coincide. However, this is not always the case.
      For example, suppose that the set $S$ has many natural automorphisms (e.g., in the case of the group $G = PGL(n,F)$,  the 
      set      ${\rm Irr}\ G$ has a~large group of automorphisms given by multiplication by characters of the group $G$).
      Then, it can happen that the correspondences $\nu$ and $\mu$ described above differ by a~twist by one of these 
      automorphisms.
      
      Now the question is ``which of the correspondences $\nu$ and $\mu$ is more correct?" This is not a~mathematical, 
      but rather a~meta-mathematical question. 

\sms      
      
      There is a~general method (widely used by Grothendieck) how to approach the  questions of this type.
      Namely, suppose that we characterized the desired correspondence by some properties. This means that we
      introduce some list $\cA$ of properties (axioms) that our correspondence between $S$ and $T$ should satisfy.
      If this system $\cA$ is good enough, then usually it is not difficult to show that the correspondence that satisfies the 
      system of axioms $\cA$ is uniquely defined -- if it exists.
      
       In this way, we can fix the desired correspondence 
      in advance, even without knowing that it exists. Then, the difficult part of the theory is to prove that such 
      correspondence does exist, but at least we know precisely what we would like to construct.
      
      \sms
      
       This is a~very powerful method Grothendieck used to produce many important constructions.
       However,  it does not completely answer the question formulated above. The reason is that we can formulate 
       two different sets of axioms $\cA$ and $\cB$, and check that they both have uniqueness property and 
       produce our correspondences $\nu$ and $\mu$. So the question ``which list of axioms is more correct?"
       is a~meta-mathematical question to  be discussed in some different (e.g., aesthetic) terms.
      
      \subsection{Basic Requirement}
      
      Langlands introduced his correspondence having in mind some analytic applications; he was motivated by and 
      mostly interested in  applications to automorphic representations and $L$-functions. 
      However, if we consider just the $p$-adic LLC, we can see that it has a~very algebraic structure. Thus, it is natural
      to formulate the properties of the desired correspondence in algebraic terms.
      
\sms      
      
      Namely, Langlands was working with representations in complex vector spaces. However, in both sides of LLC we 
      can replace the field of complex numbers by some other field. Thus, from now on, let us fix another field $\rC$
      (the field of coefficients).
      For simplicity, we assume that this field $\rC$ is  algebraically closed of characteristic $0$ --
       in other words it mimics 
      the field of complex numbers, but is a~more algebraic object, since it is not endowed with 
      any topology.
      
\sms      
      
      Now, on representation side, we consider the category $\MG$ of smooth representations
      of the group $G$ in vector spaces over the field $\rC$
      and denote by ${\rm Irr}\ G$ the set of isomorphism classes of simple objects of the category $\MG$.
      
      On the Galois side, we consider the dual group $\chG$ as a~reductive group over $\rC$, consider the ind-variety 
      ${\rm Mor}(\WF, \chG)$ of continuous (i.e., locally constant) morphisms $\phi: \WF \to \chG$ 
       and define the stack $\cLANG$ to be the quotient stack ${\rm Mor}(\WF, \chG) / \chG$.

\sms

     Notice that now both parts of the  Langlands correspondence have structure of an algebraic object over $\rC$.
      One of the main properties, that we should expect from 
      LLC  is that it is functorial in $\rC$. In particular,
      our collection of axioms for LLC should include the following \textbf{basic requirement}.
 \begin{equation*}\label{BR}
\begin{minipage}[c]{14cm}
Suppose we have two functorial constructions $S$ and $T$ that 
from a~field $\rC$ 
      produce    objects $S(\rC)$ and $T(\rC)$ in some category.
      We would like to construct a~family $\nu$ of morphisms $\nu_{\rC}: S(\rC) \to T(\rC)$.
      Then, we require that\\ \textbf{this construction commutes with all automorphisms of the field $\rC$}.
\end{minipage}
\end{equation*}

      \subsection{ }
        
        The significance  of the basic requirement is that it is \textbf{not satisfied} by the standard 
         Langlands correspondence.
        The reason is that the Langlands construction is based on the Satake isomorphism,
         and this isomorphism does not satisfy the basic requirement \ref{BR}. In other words,
         the Satake isomorphism is not completely canonical, it depends on some choices.
        We will discuss this in Subsection~\ref{wrong}.
        
\sms        
        
        The main goal of this note is to show that the basic requirement forces us to modify the stack $\cLANG$ 
        of Langlands parameters.        
        Moreover, analyzing this modified stack $\cLANG'$ we will see that there is some sign hidden in 
        the Langlands correspondence -- the sign missed in the original Langlands picture.
         
\sms         
         
          In fact, the traces of this hidden sign imbalance 
          appeared in many 
         works analyzing concrete cases of the Langlands correspondence  -- local, global and geometric (for example, see 
         \cite{DelLetter}, \cite{Clozel}, \cite{Buzzard}).  The reason is that this sign  is a~real phenomenon and ignoring 
         it makes many constructions  unnatural.
          I will discuss  this later.
          
          This work was done in the framework of the ERC grant 291612.
          It was partially done during my visits at MPIM-Bonn and IAS, Princeton. I would like to thank these 
          Institutes for very stimulating atmosphere. 
          
          I presented the results of this work in the framework of the 24-th Takagi lectures in December 2019.
          I would like to thank the Mathematical Society of Japan for inviting me to give this lecture.
          
          I would like to thank P.~Sarnak, P.~Schneider, A.~Reznikov, E.~Sayag, Y.~Sakellaridis, and S.~ Carmeli for 
          many fruitful discussions.

          \section{the Langlands description of LLC}
          
          Let me describe some arguments that lead Langlands to LLC. He already had a~construction of the dual group 
          and understood that it should play an important role in the representation theory.
         
         \subsection{{Construction of the dual group over the field $\rC$}}

Fix a~field $k$. Consider a~split torus $\bT \subset \bG$ defined over $k$. Let us denote by $X^*(\bT)$ the lattice of characters
$\lam: \bT \to G_m$ and by $X_*(T)$ the lattice of cocharacters $\nu: G_m \to \bT$, where $G_m$ is the multiplicative group 
and all morphisms are morphisms of algebraic groups. 

These two groups --- $X^*(\bT)$ and $X_*(T)$ --- are dual lattices. The value of the pairing $<\lam , \nu>$ is defined as an the integer that 
describes the composition $\lam \circ \nu: G_m \to \bT \to G_m$.
 
 The group $\bT$ can be reconstructed from any of these two lattices.
 Namely, 
 \[
 \text{$\bT = X_*(\bT) \otimes_\Z G_m$ and also
  $\bT  = {\rm Hom} (X^*(\bT), G_m)$.}
  \]
 
 \sms

Given another field $\rC$, we define the dual torus $\chT$ as a~split algebraic torus 
over $\rC$ such that
$X_*(\chT) = X^*(\bT)$, i.e., $\chT = {\rm Hom}(X_*(\tT), G_m(\rC))$.

\sms

Now, let $\bG$ be a~connected split reductive group over $k$. Let us fix a~split Cartan
subgroup $\bT \subset \bG$. Let us denote by $\Sig \subset X^*(\bT)$ the root system 
of $\bG$ and by
$\chSig \subset X_*(\bT)$ the dual coroot system of $\bG$. 
 It is the standard fact that the quadruple $ (X^*(\bT), X_*(\bT), \Sig , \chSig)$
 determines an isomorphism class of the group $\bG$.
 
\sms 
 
 Now, consider the dual quadruple $(X_*(\bT),  X^*(\bT), \chSig, \Sig)$ and construct the 
 reductive  group $\chG$ over some other field $\rC$ that corresponds to this dual quadruple. 
 By definition  this is the \textbf{Langlands dual group of $G$}.
 
 \subsection{LLC for unramified representations } \label{unramified}

The starting point of the Langlands construction is an observation that to every 
unramified irreducible representation $\pi$ of the group $G$ we can assign a~
conjugacy class $L(\pi)$ in the dual group $\chG$ 
(this is an analog of the {Artin map} in  the local class field theory).

Namely, let us fix a~hyper-special subgroup $K \subset G$ and denote by $H(G) = H_K(G)$ the corresponding 
spherical Hecke algebra with coefficients in the field $\rC$.  
 Every unramified representation $\pi$ has unique, up to scalar factor,
$K$-invariant vector $v$, and the vector $v$ defines a~character $\theta = \theta_\pi$ of the algebra $H(G)$.

According to a~Satake theorem, we can identify the algebra $H(G)$ with the algebra $O(\chT)^W$ of regular 
functions on the dual torus $\chT$ invariant  
under the action of the Weyl group $W$ of $\chG$ (it is also the Weyl group of $G$).
Thus, the character $\theta$ corresponding to $\pi$ defines a~$W$-orbit on $\chT(\rC)$, and 
hence a~semisimple conjugacy 
class of the group $\chG(\rC)$.

\sms
 
Using a~canonical morphism $p: \WF \to F^\times$ from the local class field theory, 
Langlands  interprets the above-described conjugacy class as 
an unramified morphism 
\[
L(\pi) =\phi: \WF \to \chG
\] 
defined
 up to the adjoint action of $\chG$.

\subsection{Extension to other representations -- global argument}
 
Next, Langlands uses an informal, but very powerful, global argument to suggest that
this correspondence should have a~natural  extension $L: {\rm Irr}\ G \to {\rm Lan}(G)$ to all
 irreducible representations --- we will call this $\textbf{the Langlands correspondence}$.
 
 The argument flows approximately as follows. We can consider our field $F$ as a~local counterpart of some global field $k$,
 i.e., $F = k_u$ for some place $u$ of $k$. Given an irreducible representation $\pi$ of $G$, we can hope that it can be 
 extended to a~global automorphic representation $\Pi$ of the adelic group $G(\bA_k)$.
 This representation is a~restricted tensor  product --- over all places of $k$ --- of local representations $(\pi_v, G(k_v))$.
 
 Since almost all of these representations are unramified, we can apply to them the previous construction.
  This defines a~collection of morphisms 
  \[
  \text{$\phi_v: W_v \to \chG$ , where $W_v = W_{k_v}$.}
  \]
  Now we can think about this collection as coming from one morphism $\phi: \bW \to \chG$, where $\bW$ is some 
  version of the global Galois group (we assume that the group $\bW$ contains all the Weil groups $W_{k_v}$ as subgroups).
  
  Notice that the morphism $\phi$ is uniquely determined by the collection of morphisms $\phi_v$ for almost all $v$,
  and, in its turn, it determines a~morphism $\phi_u: \WF = W_{k_u} \to \chG$ that we wanted to construct
  (this is some kind of ``analytic continuation" in the parameter $v$).
  
  I repeat: 
  this is quite informal, but very convincing argument.
  
  \subsection{What goes wrong?} \label{wrong}
  
  The correspondence described above does not satisfy the basic requirement \ref{BR}. The reason is that the
   Satake isomorphism does not satisfy this requirement. Let us analyze what goes wrong.
   
   The most transparent description of the Satake isomorphism is given by the following procedure.
   Let us fix a~Borel subgroup $B \subset G$, denote by $U$ its nil-radical and consider the natural  projection
   $p: B \to T = B / U$. This group $T$ is called \textbf{the Cartan group} of $G$. Note that it is defined by the group $G$ 
   canonically, i.e., up to unique isomorphism. 
   
    Unramified characters $\chi: T \to \rC^\times$
   naturally correspond to points of the dual torus $\chT(\rC)$. Given such a~character, we  extend it to the group $B$,
   and construct the normalized   induced representation $\pi_\chi = \mathbf{ind}_B^G(\chi)$. 
   This representation defines a~character 
   \[
   \theta_\chi: H(G) \to \rC^\times.
   \]
   Thus, an element $h \in H(G)$ defines a~function on points $\chi \in \chT$ of the dual torus, and the Satake theorem 
   is the statement that this correspondence gives an isomorphism of $H(G)$ with the algebra $O(\chT)^W$.   
   
   The problem is that the normalized induction is not compatible with the basic requirement \ref{BR}. Indeed, the usual induction is 
   clearly compatible with the basic requirement \ref{BR}, but the normalized induction differs from it by a~twist by a~character $\delta$ 
which is a~square root of the modulus character. The choice of this square root is not canonical. It is easy to make 
   it canonical if we fix a~square root $\sqrt{p} \in \rC$ of the prime number $p$ equal to the residual characteristic of $F$.
    However, this ``canonical'' choice is not 
   invariant with respect to automorphisms of the field $\rC$.
   
       Thus, we see that \textbf{the parabolic induction and the Satake isomorphism are not compatible with the basic requirement} \ref{BR}.
       
       Note that we can not simply replace normalized induction by the usual induction since having done so 
we will not get 
       $W$-invariant functions on $\chT$.

\section{Modification of the Langlands correspondence}

In this section I will describe a~modification $L'$   of the Langlands correspondence which is  compatible with the basic requirement
 \ref{BR}. Following the idea by Langlands, I will first define the correspondence $L'$ on a~large collection of 
 irreducible representations $\pi$.  This large collection consists of  representations of ``the principal series types" ---
 it contains all unramified representations. Then, I conjecture that it should naturally extend to all irreducible representations.   

\subsection{The Jacqet functor }  \label{Jacquet}
 Fix a~point  of the flag variety $X$ and denote by $B$ its stabilizer in $G$.
We denote by $U$ the unipotent radical of $B$ and by $T$ the quotient toric group $T = B/U$.

We will use the \textbf{non-normalized} Jacquet functor 
\[
\text{$J: \MG \to \MT$ given by $(\pi,V)  \mapsto (J(\pi), J(V))$,}
\]
where $J(V):= V_U$ is
 the quotient space of coinvariants, equipped with the natural action of the group $T$. 
 
  We say that an irreducible representation $(\pi,V)$ is of \textbf{the principal series type} 
  if $J(\pi)  \neq 0$.
  
\sms  
  
  Let $(\pi, V)$  be an irreducible representation of this type. Then, the space $J(V)$ is finite-dimensional.
  Let us denote by $\mu = \mu(\pi)$  
  the finite collection of all characters $\mu_i$ of the abelian group $T$ that appear in the 
  decomposition of this space.  
  
  Since $T$ is split, i.e., isomorphic to a~product of copies of the group $F^\times$,  we can interpret any
  morphism $\chi: T \to \rC^\times$ as a~morphism $\chi: F^\times \to \chT(\rC)$. 
  The Local Class Field Theory implies that we can interpret $\chi$ as a~morphism $\chi: \WF \to \chT(\rC)$.

  Thus, starting form the representation $(\pi, V)$, we have constructed a~family of morphisms $\mu_i: \WF \to \chT(\rC)$.
  Unfortunately, these morphisms are not conjugated by the action of the Weyl group W on $\chT$,
   so we can not construct from
  them one morphism $ \WF \to \chG$ defined up to conjugation. 
  
  \subsection{The dot action and normalization} \label{dot}
  
  In order to construct a~morphism $\WF \to \chG$, we will use the fact that the characters $\mu_i$ are related by a~symmetry  called the \textbf{dot action}. Let us describe it.
   
   First of all, given any (algebraic) weight $\nu: \bT \to G_m$, we can construct a~character 
   $|\nu |: T \to \rC^\times$ by composing 
   with the standard modulus morphism $| - | : F^\times \to \rC^\times$.
   Now, define the dot action of the Weyl group $W$ on the group $\text{Char} (T)$ of characters of $T$ by the condition that for every
   simple root $\al$, the dot action of the corresponding simple reflection $\sig_\al$ is given by
\[
\sig_\al \cdot (\chi) := \sig_\al(\chi) / |\al | . 
\]
   
   The following claim is the standard consequence of the theory of intertwining operators.
   
   \begin{claim} All the characters $\mu_i$ of the group $T$ lie on the same orbit of the dot action of the Weyl group $W$.
   \end{claim}
   
     This claim is just a~reformulation of the fact that  representations obtained by normalized induction from 
     two characters have a~common sub-quotient if and only if these characters  lie on the same $W$-orbit.

   \subsection{Weights of type $\rho$. } \label{rho}
   
   Let $r \in X^*(\bT)$ be a~weight. We say that $r$ is a~\textbf{weight of type $\rho$} if for any simple root $\al$ 
   we have $\sig_\al(r) = r / \al$ (here we write the operation in $X^*(\bT)$ multiplicatively).
   
   For example, if the group $G$ is semi-simple and simply connected, then the weight $\rho$ equal to the half sum 
   of positive roots has this property.
   
  It is easy to see that weights of type $\rho$  are obtained from one another via multiplication
   by algebraic characters of the group $\bG$.
   
   Suppose  we fixed a~weight $r$ of type $\rho$. Given an irreducible representation $\pi \in {\rm Irr}\ G$ of the principle series type,
    we can construct  a~collection  of characters $\mu = \mu(\pi)$ of the group $T$ as in Subsection~\ref{Jacquet}. 
    Then, we construct a~new collection $\lam = \lam(\pi)$ of characters $\lam_i$ by setting $\lam_i := |r| \cdot \mu_i $.
    
    The claim \ref{dot} means that all the characters $\lam_i$ are conjugated under the action of $W$.
    This means that the corresponding morphisms $\lam_i: \WF \to \chT(\rC)$ are conjugated under the action of $W$.
    Hence,  up to a~conjugation by 
    the group $\chG(\rC)$, they define a~morphism $\lam: \WF \to \chG(\rC)$ --- that was our goal.
    
    \subsection{Construction of the correspondence $L'$ for representations of principal series type}
    
       In general, our  group $\bG$ does not have a weight of type $\rho$. But we can remedy this.
       
       Suppose we constructed a~central extension of the algebraic group  $\bG$
       \[
       1 \to  G_m \to \bR \to \bG \to 1
       \]
       Then, by Hilbert's Theorem 90 we have a~central extension of the topological group $G$ 
       \[
       1 \to F^\times \to R \to G \to 1.
       \]
       
       Thus, we can describe any representation $\pi$ of $G$  as a~representation of the group $R$ 
       trivial on $F^\times$.
       
       Now I claim that  any group $\bG$ has a canonical central extension of this type
       and this extension 
   has a~canonical weight $r$ of type $\rho$.  We will denote this extended group by $\tbG$, the group of its 
       $F$-points by $\tG$, and its dual group over $\rC$ by $\GD$.
       
       Given an irreducible representation $\pi \in {\rm Irr}\ G$, we will lift it to a~representation of the group $\tG$; 
       then, using the canonical weight $r$, we construct a~morphism $\lam: \WF \to \GD$ defined 
       up to a~conjugation by the group  $\chG$. This is the idea of our construction of the 
       \textbf{modified Langlands parameter} $\psi = \psi(\pi): \WF \to \GD$.

 \subsection{The geometric construction of the extended group $\tbG$.}

 Let us first give a~geometric description of the group $\tbG$. We denote by $\bX$ the flag variety of 
 the group $\bG$ (we consider $\bX$ as a~$k$-scheme; its $k$-points are Borel subgroups of $G$). 
 The group $G$ acts on the set  $X = \bX(k)$ transitively.
 
\subsubsection{Square roots} 
 
 Let $\bDel$ be the anti-canonical line bundle on $\bX$; the group $\bG$ naturally acts on $\bDel$.

 \begin{lemma}  \label{sr}
 \textup{1)}  The bundle $\bDel$ has a~square root, i.e., there exists a~line bundle
 $\bL$ on the $k$-scheme $\bX$ such that its square $\bL^{\otimes 2} := \bL \otimes \bL$ is 
 isomorphic to $\bDel$.  
 
 \textup{2)}   Any two square roots $\bL$ and $\bL'$ are isomorphic as line bundles.
 
 \textup{3)}   For any line bundle $\bL$, we have ${\rm End}(\bL) = k$. In particular, the group
 of automorphisms of $\bL$ is isomorphic to $k^\times$.

 \end{lemma}
 
 \begin{remark}\normalfont Two square roots $\bL$ and $\bL'$ are not always isomorphic as square roots,
  (i.e., there is no isomorphism between these vector bundles compatible with isomorphisms of their squares with $\bDel$).
    In fact, it is easy to see that the group $k^\times$ transitively acts on the set of isomorphism classes of square roots; the kernel of this action coincides with $(k^\times)^2$, i.e., the set of isomorphism classes of square roots 
    has the size of the set $k^\times / (k^\times)^2$.
    
    \end{remark}
    
    In this lemma we consider all vector bundles and their morphisms over the field $k$. When $k$ is algebraically closed, 
    these results are standard. The general case  is an easy exercise in Galois cohomology.  

\sms

\subsubsection{ }

 Let us fix a~square root $\bL$ of the bundle $\bDel$, and use it to construct the extended group $\tbG$.
Let me describe the set of its points over the algebraic closure of $k$.
 
 Notice that for  any element  $g \in \bG$, the line bundle  $g^*(\bL)$ is isomorphic to $\bL$ since its 
 square is isomorphic to $g^*(\bDel) \simeq \bDel$.
 
\sms
 
 Now define the group $\tbG$ as a~set of pairs $(g,\al)$, where $g \in \bG$ and $\al$  is an isomorphism
 $\al: \bL \to g^*(\bL)$. Then, $\tbG$ is a~group under the natural composition law. Clearly, $\tbG$ has a~central
  subgroup $G_m$ and we have an exact sequence of algebraic groups 
 \[
 1 \to  G_m \xrightarrow{i'} \tbG \xrightarrow{p'} \bG \to 1.
 \]
 
 In particular, we see that the group $\tbG$ is connected.
 
 The square root $\bL$ defines a weight $r$ of the Borel group -- this is the canonical weigh $r$ that we wanted to construct
 
It is clear that the group $\tbG$,  all morphisms, and the weight $r$ are defined over the field $k$. 
 From Lemma ~\ref{sr} it is obvious 
 that the groups $\tbG$ constructed from two different square
 roots are \textbf{canonically} 
 isomorphic, so this construction does not depend on the choice of the square root.
 
 \subsubsection{Combinatorial description of the extended group $\tbG$}

 Let $(X^*(\bT), X_*(\bT), \Sig , \Sig\spcheck)$ be the quadruple that defines the group $\bG$.
 Consider the lattice $X^*(\tbG) := X^*(\bG) \oplus \Z r$ and extend the action of the Weyl group $W$ to this lattice
 by the condition that for every simple root $\al$, we have 
 \[
 \sig_\al(r) = r / \al.
 \]
 
 This defines all the required combinatorial data for the group $\tbG$.

 \subsubsection{Features of the extended group $\tbG$ } \label{features} {}~{}

$\bullet$ We have a~canonical central extension 
 \[
 1 \to  G_m \xrightarrow{i} \tbG \xrightarrow{p} \bG \to 1
 \] 
 and the corresponding central extension on the level of $F$-points 
 \[
 1 \to F^\times \to \tG \to G \to 1.
 \]
 
$\bullet$  We have a canonical weight $r \in X^*(\tbG)$. The composition $r \circ i : G_m \to G_m$ is identity.
 
$\bullet$   We also have a~canonical morphism $j: \tbG \to G_m$ of algebraic groups over $F$
  such that the composition 
  $j \circ i : G_m \to G_m$ is given by squaring map $z \mapsto z^2$.

Let us describe this morphism $j$ over an algebraically closed field.
  
   We fix an isomorphism of $\bL^{\otimes 2}$ with $\bDel$.  Then,  every element $g \in \tbG$ 
  induces an isomorphism $\bDel \to g^*(\bDel)$. This isomorphism differs from the geometric 
  action of the element $p(g)$ on the bundle $\bDel$ by some constant factor $j(g) \in k^\times$, and this defines the morphism
  $j$.
   
 It is clear that the morphism $j$ does not depend on the choice of the above isomorphism and is defined over the base field $k$.

Let us denote by $\bH$ the kernel of the morphism $j$. This is a~central two-fold cover of the group $\bG$ --
it might be connected or disconnected.
  
  \begin{remark}\normalfont
  The group $\bH$ is what we really are interested in. However, since this group might be disconnected,
  it is not clear how to define its dual group. So our construction in some sense is a~way to define this dual group 
  using the resolution 
  \[
  1 \to \bH \to \tbG \to G_m \to 1.
  \]
  \end{remark}
  
  \subsubsection{Langlands dual data} \label{data}
 
 Now we can consider the dual objects of all these morphisms.
 
 Denote by  $\cGD$  the dual group of the group $\tbG$, and consider the exact sequence of algebraic groups over $\rC$, 
which is dual to the sequence in Subsection~\ref{features}, namely the sequence
\[
1 \to \chG \xrightarrow{i} \cGD  \xrightarrow{p} G_m \to 1. 
\]
 
 We denote by $r$ the coweight of the group $\cGD$ corresponding to the weight $r$ above.
 We also  consider  the central morphism $j: G_m \to \cGD$
 dual to the morphism $j: \tbG \to G_m$ from Subsection~\ref{features}.
 We see that the composition $p \circ r : G_m \to G_m$ is the identity morphism,  and the 
 composition $p \circ j : G_m \to G_m$ is given by squaring.
 
 \begin{definition}\normalfont We call the collection $ (\cGD, i, p , j, r)$  the \textbf{Langlands dual data} for the group $G$.
  \end{definition}
  
  \begin{remark}\normalfont For $\bG = PGL(2)$, the group $\tbG$ is isomorphic to $GL(2)$. 
     
      The Langlands data is 
      \[
      \text{$\chG = SL(2,\rC)$, $\GD = GL(2, \rC)$, $p = \det$},
      \]
       and  the morphism 
     $j: \rC^\times \to \GD = GL(2, \rC)$ is the usual central morphism.
     
      In standard coordinates, the coweight
     ${r \in X_*(\GD)}$  is given by $(1,0)$.
     
     \sms
     
   It is very useful to play with all the constructions I described 
    for this particular case since it contains many of the subtleties of the general case.
    \end{remark} 
 
 \subsection{Modified Langlands parameters} \label{modified}
 We define the modulus character
 $\om: \WF \to \rC^\times$ as a~composition of the canonical projection $p: \WF \to F^\times$ from the local class field 
 theory and the modulus character  $ F^\times \to \rC^\times$ given by $a \mapsto |a|$. 
 In other words, $\om$ is an unramified character of 
 the Weil group $\WF$ such that  $\om(Frob) = q$, where $Frob$ is the Frobenius element and $q$ is the 
 cardinality of the residue field of $F$.
 
 Tracing our constructions it is easy to check that to every irreducible representation $\pi \in \MG$ of principal series type 
 we have assigned a~morphism
 $\psi: \WF \to \cGD$, defined up to a~conjugation by $\chG$, that satisfies the following condition:
\begin{equation*}\label{omega1} 
 \begin{minipage}[c]{14cm}
$(\omega)$ The composition of $\psi$ with the projection $p: \GD \to \rC^\times$\\
\phantom{XX}  coincides with the modulus morphism $\om$.
\end{minipage}
\end{equation*}

 Let us consider the variety ${\rm Mor}_\om(\WF, \GD)$ of all continuous  morphisms 
 $\psi: \WF \to \GD$ satisfying condition $(\om)$; define
  the \textbf{modified Langlands set} $\LANG'$ as the quotient set ${\rm Mor}_\om(\WF, \GD) / \chG$, and the 
  \textbf{modified Langlands stack} $\cLANG'$ as the quotient stack ${\rm Mor}_\om(\WF, \GD) / \chG$.

  \subsubsection{ }
  
  Then, we can formulate the modified local Langlands correspondence as follows:
 \begin{equation*}\label{omega} 
 \begin{minipage}[c]{14cm}
(\textbf{LLC}$'$) \ \  There should be a~canonical correspondence between the set ${\rm Irr}\ G$\\ 
\phantom{XXXXX} and the 
set of modified Langlands parameters
  $\psi: \WF \to \GD$ given\\ 
  \phantom{XXXXX} by points of the modified Langlands stack.
\end{minipage}
\end{equation*} 
  
  
     The advantage of this formulation is that it is compatible with the basic requirement \ref{BR}.
       I believe that the correspondence $L'$ is a~more ``correct" version of the original Langlands correspondence $L$.
   In what follows I will try to describe some consequences of this new point of view.
     
     \subsubsection{Example}
     
     Consider the group  $G = PGL(2, F)$. In this case,  the extended group is   
     $\tbG = GL(2)$ and $\GD = GL(2, \rC)$.
     
     Thus, the modified Langlands parameters $\psi$ are morphisms $\psi: \WF \to GL(2, \rC)$ that satisfy the condition
\[
\text{    ($\om$)\ \ \ \ \ \ \ \ \ \ \ \ $\det(\psi(u)) \equiv \om(u)$  for any $u \in \WF$,}
\]
defined up to a~conjugation by the group $\chG = SL(2, \rC)$.

    \sms
    
    Notice that the set $\LANG'$ of modified Langlands parameters is very close to the set $\LANG$
    of Langlands parameters,
    that is the set of morphisms $\phi: \WF \to GL(2, \rC)$ satisfying the condition  
      $\det(\phi(u)) \equiv 1$.

\sms

    If we choose a~character $\del: \WF \to \rC^\times$ such that $\del^2 = \om$, then we can identify these 
  two sets of parameters by $\psi = \del \cdot \phi$.
    
    The difference is that the modified Langlands correspondence satisfies the basic requirement \ref{BR}, while the original 
    Langlands correspondence does not have this property, since the choice of the character $\del$ is not 
    preserved by automorphisms of the field $\rC$.

\section{The Langlands dual group as an enhanced group} \label{enhanced}

\subsection{Description of the Langlands-data group $\GD$}

   Let us study in more detail the structure of the Langlands data described in Subsection~\ref{data}.
   We have a~morphism of algebraic groups $p: \GD \to G_m$; 
   its kernel is the Langlands dual group $\chG$.
   We also have a~morphism $j: G_m \to \GD$ such that the composed morphism $p \circ j: G_m \to G_m$ 
   is given by squaring.
   
   From this we immediately deduce the following
   
   \begin{proposition}
   
   The group $\GD$  is generated by the subgroups $\chG$ and $j(G_m)$. The natural morphism 
   $G_m \times \chG \to \GD$ is a~central two-fold covering, i.e., the group $\GD$ is a~quotient
  of the group $G_m \times \chG$ by a~subgroup of order $2$ generated by a central element
  ${(-1, \eps) \in G_m \times \chG}$.
  
  \end{proposition}
  
  \subsubsection{The special element $\eps \in \chG$}
  
    The implication is that the dual group has an additional structure -- it has a~distinguished central
     element $\eps$ of order $1$ or $2$. In other words, it is equipped with a~canonical morphism 
     $\eps$ of the cyclic group $\mu_2$ to the center of $\chG$, where $\mu_2 \subset C^\times$
     is the group of square roots of $1$.
     
     \begin{definition}\normalfont An \textbf{enhanced group} $(D, \eps)$ is a group $D$ 
     enhanced with a~morphism
      ${\eps: \mu_2 \to Z(D)}$,        where $Z(D)$ is the center of $D$.
      
      \end{definition}
      
      Thus, we see that the dual group has a~\textbf{canonical} structure of an enhanced group.

     This enhanced structure plays a~crucial  role in my exposition. It can be described in many different ways.
     Here is one of the descriptions.
     
      Consider the weight $t \in X^*(\bG)$  equal to the sum of all positive roots of $\bG$. 
      We can interpret $t$ as an  element of $X_*(\chT)$, i.e., as a~morphism 
     $t: G_m \to \chT$. Let $\eps$ be the restriction of this morphism to the subgroup 
     $\mu_2 \subset G_m$. 
       
        It is easy to check that the morphism  $\eps$ is 
     invariant under the action of the Weil group of $\chG$, and hence
      its image lies in the center of $\chG$.
     
     This morphism $\eps: \mu_2 \to \chG$ has been described in \cite{DelLetter}.

     \subsection{A Remark on Geometric Satake Theorem}
     
     When I discussed this additional structure with other mathematicians, they objected
     that in the Geometric Satake Theorem one constructs the dual group $\chG$ without additional structures.
     It took me some time to understand that in the Geometric Satake Theorem, the group $\chG$ also is an enhanced group.
     Unfortunately, in most of papers that I read about this, it is very difficult to discern this structure from the 
     statements of results; 
     the structure only becomes apparent when you dive into the proofs.
     
   \subsubsection{ } The standard way to formulate the Geometric Satake Theorem runs as follows. 
   Starting from the group $G$, we   construct the affine Grassmannian $\BG$ and some category 
   $\cP$ of perverse sheaves on this Grassmannian.
   
   In \cite{L}, Lusztig has defined a~convolution operation $*$ on the category $\cP$; 
   this operation endows  $\cP$ with    the structure of a~monoidal category.
   Later,  Drinfeld used the global version of the affine Grassmannian to reinterpret  the operation $*$ 
   as a~fusion product. This implies that, in fact, the category  $\cP$ has a~canonical tensor stricture; 
   i.e., a~ monoidal structure and a~symmetry constraint. 
     
   V. Ginzburg in \cite{Gi},  I. Mirkovich and K. Vilonen in \cite{MV}, used Tanakian formalism to relate 
     this category to the dual group $\chG$.
        
        The idea is as follows. Consider a~fiber functor $F: \cP \to \text{\text{Vect}}$ defined by 
        \[
        P \mapsto F(P):= H^*(\BG, P).
        \]
       Tanakian formalism implies that the tensor  category $\cP$ is equivalent to the category 
      ${\rm Rep}(R)$ of representations of some group $R$. One  computes  this 
     group and sees that it is isomorphic to $\chG$.
     
     \subsubsection{The actual construction}
     
     The above arguments are misleading. The actual situation is slightly more sophisticated. 
   Namely, the functor $F$ is not  a~fiber functor with respect 
     to the natural Drinfeld tensor structure.  Indeed, it preserves the monoidal structure,
     but is not compatible with the symmetry constraints.  In retrospect, this should be obvious,
      since cohomologies are not spaces, but \textbf{super} spaces.
     
     The  functor $F$ is, in fact, a~\textit{fiber functor} if we consider it as a~functor 
     $F: \cP \to {\rm SVect} $ into the category of super vector  spaces.
     
     Deciphering the proofs, one can see that the result  proven in \cite{MV} can be formulated as follows.
     
     \begin{proposition} \textup{1)} Consider the enhanced group $(\chG, \eps)$. Consider the category of representations 
     $(\sig, \chG, V)$  of this group in the  category ${\rm SVect}$ of finite-dimensional super vector spaces 
     that satisfy the following condition
   \[
   \textbf{The operator $\sig(\eps)$ defines the  parity on the super space $V$.}
   \]
  
   These representations form a~symmetric tensor category $\cM$ and the forgetful functor
      $F': \cM \to {\rm SVect}$ is a~fiber functor.
      
   \textup{2)}  The category $\cM$ with the fiber functor $F'$ is canonically equivalent to the category $\cP$ with the  fiber functor $F$.
      
      
      \end{proposition}
      
      Thus, we see that in the Geometric Satake Theorem, the dual group has the natural structure of an enhanced group.
      In fact, if one looks more carefully, one discerns  in this theory also the full Langlands data -- the group $\GD$ and 
      morphisms $i, p, j, r $.

     \subsection{Another description of modified Langlands parameters}

    We have seen that the Langlands-data  group $\cGD$ can be described as the quotient $\GD = G_m \times \chG /(-1, \eps)$.
    Hence, we have a~canonical covering morphism of algebraic groups $ G_m \times \chG \to \GD$.
    We can use it to give another description of the modified Langlands parameters.
    
    Namely, the squaring map defines a~two-fold central cover on the group $\rC^\times$. Using the modulus morphism
     $\om: \WF \to \rC^\times$ we get an induced covering $\tWF \to \WF$. We denote by $\eps_W$ the central element in
     $\tWF$ that generates the kernel of this covering.
     
     Now it is clear that a~modified Langlands parameter $\psi: \WF \to \GD$ determines, and is completely determined 
     by, the morphism $\psi: \tWF \to \chG$ that is \textit{genuine}, i.e., it maps $\eps_W$ to the distinguished central element
     $\eps$ of the Langlands dual group $\chG$.  
     
     Thus, we can define the set  $\LANG'$ of modified 
     Langlands parameters to be the set ${\rm Mor}_g(\tWF, \chG)$ up to the adjoint action of the group $\chG$,
     where ${\rm Mor}_g$ stands for the set of genuine morphisms. 
     Similarly, we can describe the modified Langlands stack $\cLANG'$.
     
     \begin{remark} \normalfont In fact, the covering $\tWF \to \WF$ splits, i.e., $\tWF$ is isomorphic to $\mu_2 \times \WF$.
     If we fix this splitting, we can identify genuine morphisms $\tWF \to \chG$ with all morphisms $\WF \to \chG$, i.e., identify 
     modified Langlands parameters with the usual ones.  
     
     However,  there is no canonical splitting, so this identification is not compatible with the basic requirement $\ref{BR}$.
     In fact, the choice of the positive square root of $p (=$ residual characteristic of $F$)
     defines a splitting of this covering that gives the usual Langland correspondence.
     \end{remark}
     
     \begin{remark}\normalfont
     From the Local Class Field Theory we know that there exists a~unique non-trivial two-fold central extension
     $p: W_2 \to \WF$.  It can be described as an extension induced from the double covering 
     $\Fbar^\times \to \Fbar^\times$ by the canonical morphism $p: \WF \to F^\times \to \Fbar^\times$.
     
\sms     
     
     I believe that this  extension $W_2$ should play an important role in the Langlands correspondence,
      but so far I was not able to bring it to the picture.
      \end{remark}

   \section{The absolute Satake isomorphism}
   
   Using the Langlands dual data described above we can give a~more ``correct" description of the Satake isomorphism,
   namely a~description compatible with the basic requirement~\ref{BR}.
   
   Let $G$ be a~split reductive group over a~$p$-adic field $F$, let $\rC$ be a~field of coefficients described 
    above; let    $H(G) = H_K(G)$ be
   the spherical Hecke algebra of $G$ with coefficients in $\rC$, see Subsection~\ref{unramified}.
    
   Consider the dual Langlands data $p: \GD \to G_m$.
   
   Let $q$ be the cardinality of the residue field of $F$. We consider $q$ as an element of $G_m$, and consider 
   the subvariety ${X_q = p^{-1}(q) \subset \GD}$. The dual group $\chG$ acts on this variety on the left, on the right, and by the adjoint action.
   Let me consider the adjoint action.
   
   The discussion above implies the following
   
   \begin{theorem} \textup{(\textbf{The absolute Satake isomorphism})}
   There exists a~canonical isomorphism of the Hecke algebra $H(G)$  with the algebra $O(X_q)^{\chG}$ 
   of regular functions on the variety $X_q$ invariant under the adjoint action of the group $\chG$.
   \end{theorem}
   
   Note that the variety $X_q$ together with the actions of the group $\chG$ is isomorphic to the 
   group $\chG$ which can be thought  of as the variety $X_1 = p^{-1}(1)$.
   However, this isomorphism is not canonical. In particular, if we consider 
   a~field $\rC$ that is not algebraically closed, then the absolute Satake isomorphism holds while 
   the usual Satake isomorphism does not work.
   
   \subsection{The absolute Satake isomorphism for the universal Hecke algebra}
   
  Denote by $\Lam^+$ the set of dominant coweights of $G$. The Cartan decomposition for the group $G$ states that
  $G$ is a~disjoint union of double cosets $G = \coprod  B_\lam$ under the left and right actions 
  of the compact group $K$. 
  
  Let us denote by $e_\lam \in H(G)$ the normalized bi-$K$-invariant measure on the open subset $B_\lam$. 
  These elements form a~basis of the linear space $H(G)$ and multiplication in this basis is given by a~ collection 
  of coefficients
  $a^\nu_{\lam \mu}$ via $e_\lam \cdot e_\mu = \sum a^\nu_{\lam \mu} e_\nu$.
  
  It is a~standard fact that these coefficients are polynomials in $q$. This means that there exists a~family of polynomials 
  $P^\nu_{\lam \mu} \in \Z[x]$, 
  defined purely in terms of the quadruple $(\bX^*(G),\bX_*(G), \Sig, \Sig\spcheck)$, such that 
  $a^\nu_{\lam \mu} \equiv  P^\nu_{\lam \mu}(q)$.
  
  Hence, we can consider the universal Hecke algebra $\cH$ over the algebra  ${A = \rC(q,q^{-1})}$ 
  of Laurent polynomials 
  with the basis $e_\lam$ and multiplication defined by the collection of polynomial $P^\nu_{\lam \mu}$.
  
  The absolute Satake isomorphism can be reformulated as follows.
  
  \begin{theorem} There exists a~canonical isomorphism of the universal Hecke algebra $\cH$ with the algebra 
  of regular functions on the group $\GD$ invariant under the adjoint action of the group $\chG$. The structure of the 
  $A$-algebra is given by the natural morphism 
  \[
  {p^*: A = O(G_m) \to O(\GD)}.
  \]
  
 \end{theorem}
 
 In fact, this isomorphism can be defined over $\Q$. Probably, after taking care of some normalizations, 
 it would be defined over $\Z$ as well.

\section{Functoriality and almost algebraic groups}

\subsection{The Langlands Functoriality }
Consider two reductive groups $H$ and $G$; suppose that we are given a~
morphism of algebraic groups $\nu: \chH \to \chG$. Then, any Langlands parameter $\phi: \WF \to\chH$ defines 
a Langlands parameter $\phi \circ \nu : \WF \to \chG$. Since Langlands parameters are related to representations, 
Langlands formulated the following functoriality principle.

\begin{principle} Given a~morphism of algebraic groups $\nu: \chH \to \chG$, there should exist a~correspondence
 $\nu_*$ from the set $\IrrH$ to the set  ${\rm Irr}\ G$.
The same is true for automorphic representations.
\end{principle}

This principle is a~very powerful tool since in many cases one can try to establish such correspondence directly
without trying to understand the mystery of the dual groups.

Now consider this matter from the perspective of the modified Langlands correspondence.
Both $\chH$ and $\chG$ are enhanced groups. If 
$\nu$ is a~morphism of enhanced groups, then
the functoriality is easy to accept. But suppose that $\nu$ is not compatible with enhancements.
 Then, it is not clear what to do.
 
 \begin{example} \normalfont Let $G = PGL(2,F)$, $H = T$ its split torus. We have the standard embedding 
\[
\nu: \chH = \rC^\times \to \chG = SL(2, \rC). 
\]
This embedding is not compatible with the enhancement of the group $\chG$,
 so it is not clear that we should  have a~nice correspondence $\nu_*$ from $\IrrH$ to ${\rm Irr}\ G$.
 
 Of course, we have a~normalized induction functor $\mathbf{ind} : {\rm Irr}\ H \to {\rm Irr}\ G$,
  but we have seen that this functor is not natural since it does not satisfy the basic requirement $\ref{BR}$.
 \end{example}
 
\subsection{ } 
 
 I suggest to generalize the class of groups
 we would like to consider, so that we can formulate this principle in more
 general situations.
 
 \begin{definition} \normalfont An \textbf{almost algebraic reductive group} $(\bG, \eps)$ is a~reductive group $\bG$ 
 equipped with a~central 
 enhancement  morphism 
 $\eps: \mu_2 \to \chG$, where $\mu_2 \subset G_m$ is the group of square roots of $1$. 
 \end{definition}
 
 Given an almost  algebraic group $(\bG, \eps)$, we construct a~topological group $G'$ as follows.
 
 Consider the group $\bS = G_m \times \bG$ and the natural embedding $\eps' = i \times \eps: \mu_2 \to \chS$.
  This defines an isogeny (two-fold covering)  $p: \chS \to \chR$ of connected reductive groups  over $\rC$, 
  where $\chR:= \chS /\mu_2$. Passing to the dual groups we get a two-fold cover $p:\bR \to \bS$ of algebraic
   groups over $F$ with fiber $\mu_2$.
   
  Let us consider the algebraic group $\bH = p^{-1}(\bG)$. We have a central extension $p: \bH \to \bG$
  with fiber $\mu_2$.  Note that the group $H$ might be disconnected.

   Starting with the central extension of algebraic groups $p: \bH \to \bG$, we construct the central extension 
   of topological groups $p: G' \to G$.
 Namely, we denote by $G'$ the preimage in $\bH(\Fbar)$ of the subgroup $G = \bG(F) \subset \bG(\Fbar)$.
 
 The morphism $p: G' \to G$ is a~topological covering of topological groups, its kernel $\bA$ is canonically isomorphic 
 to the group $\mu_2$, i.e., it has a~canonical character $\xi: \bA \to \rC^\times$. 
 
 \begin{definition}\normalfont Given an almost algebraic group $(G, \eps)$, we denote by $\MGe$ the category of smooth representations of the topological group $G'$ such that the action of the subgroup $\bA$ is given by the character $\xi$
 (genuine representations of the group $G'$).
 \end{definition}
 
 Notice that this category $\MGe$ is very similar to the category $\MG$. 
 Usually, one can easily pass from one of them to the other.
 
 If $\eps$ is the standard enhancement of $G$ described in Section \ref{enhanced}, then $G' = \mu_2  \times G$, and
 hence the category $\MGe$ is canonically equivalent to the category $\MG$.
 
\sms 
 
 Now we can formulate the generalized version of the functoriality principle. 
 
 \begin{principle} \normalfont
 Consider reductive groups $\bH$ and $\bG$, and a~morphism of algebraic groups ${\nu: \chH \to\chG}$. 
 Let $\eps: \mu_2 \to \chH$ be an enhancement of the group $H$. Assume that the composition  $e\nu$ 
 defines an enhancement $e$ of the group $G$, i.e., that the image of  $\nu \circ \eps$ is central. 
 Then, there should be a~correspondence 
 between simple objects in $\MHe$ and  in $\MGe$.
 \end{principle}

 \section{Automorphic $L$-functions}
 
 One of the most important applications of the Langlands philosophy was his construction of a~family of $L$-functions
 attached to automorphic representations and representations of the dual group.
 
 I will not repeat Langlands construction, but just describe the modified construction that arises if we  
 apply the Langlands approach to the modified Langlands correspondence. In order not to confuse things, I 
 will denote these new functions by $R$, not by $L$.
 
 \subsection{Construction of partial $R$-functions}
 
 Fix a~connected reductive group $\bG$ over a~global field $k$ (like before, I assume that $\bG$ is split).
  Let us consider the adelic group $G = \bG(\bA_k)$, and its discrete subgroup $\Gam = \bG(k)$;  denote by 
  $X$ the automorphic quotient space $X = \Gam \backslash G$.
 
 Let $\cF(X)$ be the space of complex-valued functions on $X$ of moderate growth 
 (here we work with complex vector spaces, so we take $\rC = \C$).
 Fix some irreducible automorphic representation $\Pi$ of the group $G$ realized in the space $\cF(X)$.
 (The most interesting is the case of cuspidal representations.)

   It is known that the representation $\Pi$ is the restricted tensor product of local representations 
   $(\pi_v, G_v)$ over all places $v$ of the field $k$. Here, $G_v = \bG(k_v)$ is a~locally compact group
   and $\pi_v$ is its irreducible representation. Thus, starting with the automorphic representation $\Pi$, we constructed 
   a~collection of representations $\pi_v \in \text{Irr}(G_v)$. Almost all these representations are unramified. 

\sms 
 
 Now,  fix a~complex-analytic representation $\tau: \GD \to GL(E)$, 
 see Subsection~\ref{data}.
 
I would like to assign a~complex number $R(\Pi, \tau)$ to this data.
  Formally, this number is the product of local factors $R(\pi_v, \tau)$.
   
 Consider the $p$-adic field $F = k_v$.
 Using the modified Langlands 
 correspondence we can assign to the representation  $\pi_v$ a~morphism $\psi_v: \WF \to \GD$.
 Combining it with the representation $\tau: \GD \to GL(E)$ we get a~finite-dimensional representation 
 $\cR_v = \tau \circ \psi_v $ of the Weil  group $\WF$ in the space $E$. We would like to define the local factor 
 $R(\pi_v, \tau)$ in terms of this representation.
 
 Suppose $v$ is an unramified place. Then, the representation $\cR_v$
   is unramified. Therefore it is completely determined by the image of the Frobenius class
 ${X_v = \cR_v(Frob) \in GL(E)}$.   Let us postulate that in this case the local factor $R(\pi_v, \tau)$ is defined as follows
\[
R(\pi_v, \tau): = \det(1- X_v)^{-1} = \det(1 - \tau \circ \psi_v(Frob))^{-1}.
\]
 
 \sms
 
 I assume that one can define correct factors for ramified places also, but so far let us consider the partial 
 constant $R_S$,
 where $S$ is a~finite subset of places of the field $k$, a subset that contains all ramified and Archimedean places.
 Namely, we define this constant as the product of local factors  $R_S(\Pi, \tau):= \prod_{v \notin S} R(\pi_v, \tau)$.
 
 \sms
 
 This formal product is usually not convergent, but we can  regularize it in a~standard way.
 Namely,  we include each representation $\cR_v$ into a~ family of representations ${\cR_v(s) = \om^{-s} \cdot \cR_v}$,
 where $\om$ is the modulus character of the group $\WF$, see Subsection~\ref{data}. 
 
 From these representations we construct the local factors $R(\pi_v, \tau; s)$ and define 
 the partial $R$-function $$R_S(\Pi, \tau; s) := \prod_{v \notin S} R(\pi_v, \tau; s)$$.
 
 This product absolutely converges for $\Re(s) \gg 0$ and defines a~function holomorphic in~$s$.
 Following Langlands, I  conjecture that this function has  meromorphic continuation, and then 
 define 
 $R(\Pi, \tau)$  as $ R(\Pi, \tau; s)$ evaluated at $s = 0$.
 
   \subsubsection{ }
   
   Of course, this definition essentially mimics Langlands's definition of automorphic $L$-functions.
   In other words,  these $R$-functions are, probably, just the same $L$-functions with a~slightly 
   different normalization of parameters.
   However, sometimes such normalizations do matter.

      \begin{example} \normalfont For a~given $L$-function $L(s)$, it is important to study its special values. 
      In other words,
      one would like to describe special points $s$ for which   the values $L(s)$ have an arithmetic significance,
       and try to interpret these values.
      
      In the normalization  I described, it is quite clear where to look for these  spacial points. Namely, 
   if  the automorphic representation $\Pi$ is of algebraic type and the representation $\tau$ of the 
      Langlands-data group $\GD$ is algebraic then the point
      $s=0$ should be special.    In fact, in this case all integer points $s$ are special,  since we have 
      the identity
      where $p$ is the $1$-dimensional representation of  $\GD$ corresponding to the morphism 
      ${p: \GD \to \ G_m}$ in the Langlands data, see Subsection~$\ref{data}$.
      
      \end{example} 
      
      \subsubsection{Functional equation}
      
  I do not know how to describe the meromorphic continuation  of these $R$-functions
      and how to write the ramified and Archimedean factors.  However, if all this is done, then I would be able to 
      guess the shape of the functional equation.
      
      Namely, given a~cuspidal  automorphic representation $\Pi$ of $G$ and a~representation $\tau$ of the group $\GD$,
      consider the contragradient automorphic representation $\tilde{\Pi}$ and the representation 
      ${\tau\spcheck}:= p \otimes\tau^*$ of the group $\GD$. Then, the functional equation should give a~ simple relation
      between functions $R(\Pi, \tau; s)$ and $R(\tilde{\Pi}, {\tau\spcheck}; -s) $.

\printindex


\def\cprime{$'$}

    \end{document}